\documentclass{article}%
\usepackage{amsmath}
\usepackage{amsfonts}
\usepackage{amssymb}
\usepackage{graphicx}
\usepackage[frenchb]{babel}%
\begin{document}

\title{Sur l'article \textit{On some ideals of differentiable functions}}
\author{Bernard Teissier}
\maketitle

\section{Commentaires} Cet article introduit deux id\'ees importantes: la d\'efinition pour un id\'eal ${\mathcal J}$ de fonctions diff\'erentiables dans un ouvert $U$ de ${\mathbf R}^n$ d'\^etre "de \L ojasiewicz" et le fait que cette propri\'et\'e implique une condition de r\'egularit\'e pour le lieu des z\'eros des fonctions de l'id\'eal: le fait que ce lieu des z\'eros contient un ouvert dense qui est une sous-vari\'et\'e de $U$. Ce r\'esultat sera repris dans le livre \cite{To} de Tougeron.\par La seconde id\'ee est celle d'"extension jacobienne" d'un id\'eal de fonctions diff\'erentiables, dont Thom pense qu'elle "joue le r\^ole" du Nullstellensatz dans le cas r\'eel" puisqu'elle permet (Th\'eor\`eme 2 de l'article) de construire en chaque point d'un ouvert dense du lieu $E$ des z\'eros d'un id\'eal de \L ojasiewicz un id\'eal de d\'efinition de $E$. Comment ceci peut-il \^etre rapproch\'e du th\'eor\`eme des z\'eros r\'eel de Risler dans \cite{R} pour les id\'eaux de fonctions analytiques r\'eelles, qui sont de \L ojasiewicz ? \par\noindent

A la fin de l'article Thom propose l'id\'ee qu'il existe des op\'erations canoniques d'extension d'un id\'eal de fonctions analytiques g\'en\'eralisant celle d'extension jacobienne et dont les lieux de z\'eros forment une stratification de Whitney du lieu des z\'eros de l'id\'eal. Dans le cas analytique complexe, la caract\'erisation alg\'ebrique des stratifications de Whitney par l'auteur de ces lignes dans \cite{Te} est peut-\^etre une \'etape vers la r\'ealisation de ce magnifique r\^eve de Thom.\par
Enfin, on peut s'\'etonner que Thom ne cite pas l'article fondateur \cite{W} de Whitney sur les id\'eaux de fonctions diff\'erentiables. Il est vrai que l'esprit de cet article est enti\`erement diff\'erentiel alors que Thom s'int\'eresse \`a ce qui s\'epare l'analytique du diff\'erentiel, avec peut-\^etre en arri\`ere-plan la diff\'erence entre ensembles stratifi\'es diff\'erentiels et analytiques.

\end{document}